\documentclass{amsart}

\usepackage{amssymb}
\usepackage[hidelinks]{hyperref}
\usepackage[textsize=footnotesize]{todonotes}
\usepackage{tikz}
\usetikzlibrary{cd}
\usepackage{float}

\theoremstyle{plain}
\newtheorem{theorem}{Theorem}[section]
\newtheorem{lemma}[theorem]{Lemma}

\newtheorem{proposition}[theorem]{Proposition}

\theoremstyle{definition}
\newtheorem{definition}[theorem]{Definition}

\theoremstyle{remark}

\newcommand{\bN}{\mathbb{N}}
\newcommand{\N}{\bN}

\newcommand{\cA}{\mathcal{A}}
\newcommand{\cB}{\mathcal{B}}

\newcommand{\cE}{\mathcal{E}}

\newcommand{\cI}{\mathcal{I}}
\newcommand{\I}{\cI}
\newcommand{\cJ}{\mathcal{J}}
\newcommand{\J}{\cJ}

\newcommand{\cP}{\mathcal{P}}

\newcommand{\continuum}{\mathfrak{c}}
\newcommand{\cc}{\continuum}

\newcommand{\Fin}{\mathrm{Fin}}

\newcommand{\Exh}{\mathrm{Exh}} 
 
\newcommand{\lh}{\mathrm{lh}} 
\newcommand{\NF}{\mathcal{NF}}

\begin{document}

%%%%%%%%%%%%%%%%%%%%%%%%%%%%%%%%%%%%%%%%%%%%%%%%%%
%%%%% TITLE
%%%%%%%%%%%%%%%%%%%%%%%%%%%%%%%%%%%%%%%%%%%%%%%%%%

\title{On Farah ideals and beyond}

%%%%%%%%%%%%%%%%%%%%%%%%%%%%%%%%%%%%%%%%%%%%%%%%%%
%%%%% AUTHORS
%%%%%%%%%%%%%%%%%%%%%%%%%%%%%%%%%%%%%%%%%%%%%%%%%%

\author[Adam Kwela]{Adam Kwela}
\address[Adam Kwela]{Institute of Mathematics\\ Faculty of Mathematics\\ Physics and Informatics\\ University of Gda\'{n}sk\\ ul.~Wita  Stwosza 57\\ 80-308 Gda\'{n}sk\\ Poland}
\email{Adam.Kwela@ug.edu.pl}
\urladdr{https://mat.ug.edu.pl/~akwela}

\thanks{The author was supported by the Polish National Science Centre project OPUS No. 2024/53/B/ST1/02494.}

%%%%%%%%%%%%%%%%%%%%%%%%%%%%%%%%%%%%%%%%%%%%%%%%%%
%%%%% DATE
%%%%%%%%%%%%%%%%%%%%%%%%%%%%%%%%%%%%%%%%%%%%%%%%%%

\date{\today}

%%%%%%%%%%%%%%%%%%%%%%%%%%%%%%%%%%%%%%%%%%%%%%%%%%
%%%%% MSC (MATHEMATICAL SUBJECT CLASSIFICATION)
%%%%%%%%%%%%%%%%%%%%%%%%%%%%%%%%%%%%%%%%%%%%%%%%%%


%%%%%%%%%%%%%%%%%%%%%%%%%%%%%%%%%%%%%%%%%%%%%%%%%%
%%%%% KEYWORDS
%%%%%%%%%%%%%%%%%%%%%%%%%%%%%%%%%%%%%%%%%%%%%%%%%%


%%%%%%%%%%%%%%%%%%%%%%%%%%%%%%%%%%%%%%%%%%%%%%%%%%
%%%%% ABSTRACT
%%%%%%%%%%%%%%%%%%%%%%%%%%%%%%%%%%%%%%%%%%%%%%%%%%

\begin{abstract}
An ideal $\mathcal{I}$ on $\omega$ is called Farah if there are hereditary closed subsets $C_n$ of $\mathcal{P}(\omega)$, for $n\in\omega$, such that $\mathcal{I}=\{A\subseteq\omega:\ \forall_{n\in\omega}\ \exists_{k\in\omega}\ A\setminus k\in C_n\}$. Farah asked in [Trans. Amer. Math. Soc. 356 (2004), no. 6] if every $\mathbf{F_{\sigma\delta}}$ ideal is a Farah ideal. We show that this is false by providing both tall and non-tall counterexamples. Moreover, we show that there is an antichain in the sense of Kat\v{e}tov order of length $\continuum$ consisting of $\mathbf{F_{\sigma\delta}}$ non-Farah ideals.

We also investigate a closely related problem concerning sets of ideal limit points. We show that if $\mathcal{I}$ is an $\mathbf{G_{\delta\sigma\delta}}$ ideal on $\omega$ and $x\in X^\omega$ is a sequence in an uncountable Polish space $X$, then
\[
\Lambda_x(\mathcal{I})=\left\{p\in X:\ \exists_{A\notin\mathcal{I}}\ x\restriction A\text{ converges to }p\right\}
\]
is $\mathbf{F_{\sigma}}$. This answers [Top. App. 312 (2022)] (partially) and [J. Symb. Logic, DOI: 10.1017/jsl.2026.10181]. Earlier it was known to be true only in the case of Farah ideals.
\end{abstract}

%%%%%%%%%%%%%%%%%%%%%%%%%%%%%%%%%%%%%%%%%%%%%%%%%%
%%%%% MAKE TITLE
%%%%%%%%%%%%%%%%%%%%%%%%%%%%%%%%%%%%%%%%%%%%%%%%%%

\maketitle

%%%%%%%%%%%%%%%%%%%%%%%%%%%%%%%%%%%%%%%%%%%%%%%%%%
%%%%% SECTION
%%%%%%%%%%%%%%%%%%%%%%%%%%%%%%%%%%%%%%%%%%%%%%%%%%

\section{Introduction}

Before explaining main topic of this paper, we recall some rather standard terminology. 

We say that a subset of $\mathcal{P}(\omega)$ is hereditary if it is closed under taking subsets of its elements. A family $\I\subseteq\cP(\omega)$ is an ideal if it is hereditary and closed under taking finite unions of its elements. In this paper we additionally assume that all ideals are proper subsets of $\cP(\omega)$ (that is, $\omega\notin\I$) and that $\Fin\subseteq\I$, where $\Fin$ is the collection of all finite subsets of $\omega$. For an ideal $\I$ on $\omega$ we denote its co-ideal by $\I^+=\cP(\omega)\setminus\I$.

We say that ideals $\I$ and $\J$ on $\omega$ are isomorphic if there is a bijection $f:\omega\to \omega$ such that: 
\[
B\in \J \iff f[B]\in \I,
\]
for all $B\subseteq\omega$. Using the notion of isomorphism of ideals we can consider ideals on arbitrary countably infinite sets and transfer them on $\omega$ via bijections.

We treat the power set $\mathcal{P}(\omega)$ as the space $2^\omega$ by identifying subsets of $\omega$ with their characteristic functions. Thus, we can talk about descriptive complexity of subsets of $\mathcal{P}(\omega)$. In particular, all basic clopen subsets of $\mathcal{P}(\omega)$ are of the form $\{A\subseteq\omega:\ F\subseteq A,\ G\cap A=\emptyset\}$ for $F,G\in\Fin$.

An ideal $\mathcal{I}$ on $\omega$ is called Farah if there is a sequence $(C_n)_{n\in\omega}$ of hereditary closed subsets of $\mathcal{P}(\omega)$ such that 
\[
\mathcal{I}=\{A\subseteq\omega:\ \forall_{n\in\omega}\ \exists_{k\in\omega}\ A\setminus k\in C_n\}.
\]
It is immediate that each Farah ideal is $\mathbf{F_{\sigma\delta}}$. This class of ideals has been introduced in $2004$ by Farah (see \cite[Definition 2.2]{MR2048515}), who also asked whether every $\mathbf{F_{\sigma\delta}}$ ideal is a Farah ideal (see \cite[Question 15.2]{MR2048515}). Later Hru\v{s}\'{a}k and Meza-Alc\'{a}ntara provided a useful characterization:
\begin{theorem}[{\cite[Proposition 3.6 and Theorem 3.8]{MR2849045}}]
\label{char}
The following are equivalent for every ideal $\I$ on $\omega$:
\begin{itemize}
    \item[(a)] $\I$ is a Farah ideal,
    \item[(b)] $\I$ is a countable intersection of hereditary $\mathbf{F_{\sigma}}$ sets closed under finite modifications,
    \item[(c)] $\I$ is a countable intersection of $\mathbf{F_{\sigma}}$ sets closed under finite modifications.
\end{itemize}
\end{theorem}

Hru\v{s}\'{a}k and Meza-Alc\'{a}ntara in \cite[Definition 3.11]{MR2849045} introduced a closely related notion of weakly Farah ideals: an ideal on $\omega$ is weakly Farah if it is a countable intersection of hereditary $\mathbf{F_{\sigma}}$ sets. It immediately follows from Theorem \ref{char} that each Farah ideal is weakly Farah and it is obvious that each weakly Farah ideal is $\mathbf{F_{\sigma\delta}}$. The Authors of \cite{MR2849045} asked whether those implications can be reversed. In Section $2$, after proving some folklore lemmas, we show that the notions of Farah and weakly Farah ideals coincide. During final stages of preparation of this paper the same result has been independently communicated by He, Zhuang and Zhou in \cite[Theorem B]{Chiny} (with a very similar proof). 

In Section $3$ we prove our main results: we provide both tall and non-tall counterexamples to Farah's question (recall that an ideal $\I$ on $\omega$ is tall if for every $A\in[\omega]^\omega$ there is $B\in[A]^\omega\cap\I$) and construct an antichain in the sense of Kat\v{e}tov order of length $\continuum$ consisting of $\mathbf{F_{\sigma\delta}}$ non-Farah ideals (recall that $\I$ is below $\J$ in the Kat\v{e}tov order $\leq_K$ if there is an arbitrary function $f:\omega\to \omega$ such that $f[A]\in\I^+$ for all $A\in\J^+$). In particular, there are $\continuum$ many pairwise non-isomorphic $\mathbf{F_{\sigma\delta}}$ non-Farah ideals. 

Our studies are motivated by the variety of applications of Farah ideals, that appeared over recent years:
\begin{itemize}
    \item Farah ideals have been introduced by Farah in the context of gaps in quotient Boolean algebras over ideals, that is, algebras of the form $\cP(\omega)/\I$. In \cite{MR2048515} he proved, under the assumption of Proper Forcing Axiom, that all isomorphisms between quotient Boolean algebras over Farah ideals have continuous liftings. As far as we are concerned, it is expected that the above conclusion applies to all Borel ideals. However, the class of Farah ideals is still the largest, for which this has been proved.
    \item Hru\v{s}\'{a}k and Meza-Alc\'{a}ntara studied Farah ideals in the context of an infinite game involving two ideals. This game has been introduced by them in \cite{MR2849045}. It induces a pre-order $\sqsubseteq$ on ideals. They were able to show that all $\mathbf{F_{\sigma}}$ ideals are $\sqsubseteq$-equivalent and form the least equivalence class. Farah ideals are connected to the question about the number of equivalence classes of $\mathbf{F_{\sigma\delta}}$ non-$\mathbf{F_{\sigma}}$ ideals -- Hru\v{s}\'{a}k and Meza-Alc\'{a}ntara only managed to prove that there are at least two such classes and that all Farah ideals are $\sqsubseteq$-below members of the greater of those two classes.
    \item Farah ideals appeared in the context of sets of ideal limit points, that is, sets of the form:
    \[
    \Lambda_x(\I)=\left\{p\in X:\ \exists_{A\in\I^+}\ (x\restriction A\text{ converges to }p)\right\}
    \]
    where $x\in X^\omega$ is a sequence in an uncountable Polish space $X$. He, Zhang and Zhang showed in \cite[Theorem 4.5]{MR4393937} that in the case of Farah ideals the set $\Lambda_x(\mathcal{I})$ is always $\mathbf{F_{\sigma}}$ and asked in \cite[Question 4.9]{MR4393937} if this is true for all $\mathbf{F_{\sigma\delta}}$ ideals (note that \cite{MR4393937} uses a more general framework of first countable, not necessary Polish, spaces). Later it was proved in \cite[Theorem 10.7]{FKL24} that if $\I$ is a $\mathbf{G_{\delta\sigma\delta}}$ ideal, then:
    \[
    \mathcal{L}_X(\I)=\left\{\Lambda_x(\I):\ x\in X^\omega\right\}\cup\{\emptyset\}
    \]
    is equal to one of the following three: the family of all closed subsets of $X$, the family of all $\mathbf{F_{\sigma}}$ subsets of $X$, or the family of all analytic subsets of $X$. However, it was left as an open problem whether the last possibility can be realized (see \cite[Question 10.9]{FKL24}).
    \item For an ideal $\I$ on $\omega$ and a sequence of characters $v=(v_n)_{n\in\omega}$ of a Polish group $X$, consider the following subgroup:
    \[
    H_{v}(\I)=\left\{x\in X:\ \I\text{-}\lim_{n\to \infty} v_n(x)=0\right\},
    \]
    where $\I\text{-}\lim_{n\to \infty} v_n(x)=0$ means that $\{n\in\omega:\ \|v_n(x)\|\geq\varepsilon\}\in\I$ for all $\varepsilon>0$. Interestingly, every subgroup of $X$ is of this form, for suitable $\I$ and $v=(v_n)_{n\in\omega}$ (see \cite[Theorem 2.11]{filipow2026characterizedsubgroupsunitcircle}). In \cite[Theorem 2.5]{filipow2026characterizedsubgroupsunitcircle}, extending some results of \cite{protasov, MR1240629, MR4078214} and answering \cite[Problem 6.11]{MR4078214} and \cite[Problem 7.27]{MR4932230}, the Authors proved that in the case of Farah ideals the subgroup $H_{v}(\I)$ is always $\mathbf{F_{\sigma\delta}}$. It is unknown if this is true for all $\mathbf{F_{\sigma\delta}}$ ideals. 
    \item Very recently He, Zhue and Zhang proved that if $\I$ is a Farah ideal then there is a closed $\I$-MED family, that is, a closed set $\cE\subseteq\omega^\omega$ such that $\{n\in\omega:\ f(n)=g(n)\}\in\I$ for all distinct $f,g\in\cE$ and for every $h\in\omega^\omega$ one can find $f\in\cE$ such that $\{n\in\omega:\ f(n)=h(n)\}\in\I^+$ (see \cite[Corollary C]{Chiny}). In the same paper it has been shown that Farah ideals are uniformly weakly Ramsey, that is, they satisfy some Ramsey-type property, for which the witnesses can be chosen by a single Borel map (see \cite[Theorem B]{Chiny}).
\end{itemize}

In Section $4$ we prove that, despite the existence of $\mathbf{F_{\sigma\delta}}$ non-Farah ideals (established in Section $3$), it is true that $\Lambda_x(\I)$ is $\mathbf{F_{\sigma}}$ for every $\mathbf{G_{\delta\sigma\delta}}$ ideal $\I$ and every sequence $x\in X^\omega$ is an uncountable Polish space $X$. This answers \cite[Question 10.9]{FKL24} and \cite[Question 4.9]{MR4393937}. However, it should be mentioned that in the case of \cite[Question 4.9]{MR4393937} we only get a partial solution, as we work in the framework of uncountable Polish spaces (in order to be able to apply machinery developed in \cite{FKL24}). We do not know how our results affect the rest of the problems mentioned.

We use standard notation concerning $0-1$ sequences. In particular, for $s\in 2^{<\omega}\cup 2^\omega$ by $\lh(s)$ we denote its length. If $s\in 2^{<\omega}$ and $t\in 2^{<\omega}\cup 2^\omega$ then $s^\frown t$ is the concatenation of those sequences, we write $s\subseteq t$ if $\lh(s)\leq\lh(t)$ and $t\restriction \lh(s)=s$, and we write $s\perp t$ to denote the situation when $s\not\subseteq t$ and $t\not\subseteq s$. The symbols $0_n$ and $0_\infty$ stand for the sequences constantly equal to zero of length $n$ and $\omega$, respectively (in particular, $0_0=\emptyset$). For $s\in 2^{<\omega}$ we write $[s]$ to denote the set $\{x\in 2^\omega:\ s\subseteq x\}$. Note that $[s]\subseteq[t]$ if and only if $t\subseteq s$. Similar notation applies to sequences belonging to $\omega^{<\omega}$. Finally, by $\mathbb{Q}(2^\omega)$ we mean the family of all $x\in 2^\omega$ that are constantly equal to zero from some point on. 

%%%%%%%%%%%%%%%%%%%%%%%%%%%%%%%%%%%%%%%%%%%%%%%%%%
%%%%% SECTION
%%%%%%%%%%%%%%%%%%%%%%%%%%%%%%%%%%%%%%%%%%%%%%%%%%

\section{Farah ideals and weakly Farah ideals are the same}

\begin{proposition}[Folklore]
\label{prop:hered_closed}
Every non-empty hereditary closed subset of $\mathcal{P}(\omega)$ is of the form:
\[
\left\{A\subseteq\omega:\ \forall_{F\in\cB}\ F\not\subseteq A\right\},
\]
for some family $\cB$ of non-empty finite sets.
\end{proposition}

\begin{proof}
Let $C$ be a non-empty hereditary closed subset of $\mathcal{P}(\omega)$. If $C=\cP(\omega)$ then $\cB=\emptyset$ works. Thus, we will assume that $C\neq\cP(\omega)$. There are families of finite sets $\{F_n:\ n\in\omega\}$ and $\{G_n:\ n\in\omega\}$ such that:
\[
C=\bigcap_{n\in\omega}\left\{A\subseteq\omega:\ F_n\not\subseteq A\ \vee\ G_n\cap A\neq\emptyset\right\}.
\]
Since $F_n\cap G_n\neq\emptyset$ implies $\left\{A\subseteq\omega:\ F_n\not\subseteq A\ \vee\ G_n\cap A\neq\emptyset\right\}=\cP(\omega)$, we can assume that $F_n\cap G_n=\emptyset$ for all $n$.

We will show that $C=\left\{A\subseteq\omega:\ \forall_{n\in\omega}\ F_n\not\subseteq A\right\}$. The inclusion $\supseteq$ is obvious. For $\subseteq$, suppose that $A\in C$ and fix $n\in \omega$. Assume towards contradiction that $F_n\subseteq A$. Then $G_n\cap A\neq\emptyset$. Since $C$ is hereditary, $A\setminus G_n\in C$, but $G_n\cap (A\setminus G_n)=\emptyset$ and $F_n\subseteq A\setminus G_n$ (by $F_n\cap G_n=\emptyset$), a contradiction.

Since $F_n=\emptyset$ would imply $C\subseteq\emptyset$, we conclude that $\cB=\{F_n:\ n\in \omega\}$ works.
\end{proof}

\begin{proposition}[Folklore]
\label{prop:hered_F_sigma}
Every non-empty hereditary $\mathbf{F_\sigma}$ subset of $\mathcal{P}(\omega)$ is of the form:
\[
\left\{A\subseteq\omega:\ \exists_{n\in\omega}\ \forall_{F\in\cB_n}\ F\not\subseteq A\right\},
\]
for some sequence $\left(\cB_n\right)_{n\in\omega}$ of families of non-empty finite sets such that $\cB_n\supseteq\cB_{n+1}$ for all $n\in\omega$.
\end{proposition}

\begin{proof}
Let $H\subseteq \mathcal{P}(\omega)$ be a non-empty hereditary $\mathbf{F_\sigma}$ set. Then $H=\bigcup_{n\in\omega}L_n$, where each $L_n$ is compact, non-empty, and $L_n\subseteq L_{n+1}$ for all $n$. 

Define $\phi:\cP(\omega)\times\cP(\omega)\to\cP(\omega)$ by $\phi(A,B)=A\cap B$ and let $M_n=\phi[L_n\times\cP(\omega)]$ for all $n$. Then each $M_n$ is compact (by continuity of $\phi$). Observe that $M_n=\left\{A\subseteq\omega:\ \exists_{B\in L_n}\ A\subseteq B\right\}$. Thus, $(M_n)_{n\in\omega}$ is an $\subseteq$-increasing sequence of non-empty hereditary sets. Moreover, $H=\bigcup_{n\in\omega}M_n$, since $H=\bigcup_{n\in\omega}L_n\subseteq\bigcup_{n\in\omega}M_n\subseteq H$ (by heredity if $H$).

By Proposition \ref{prop:hered_closed}, for each $n\in\omega$ there is a family $\cB'_n$ of non-empty finite sets such that $M_n=\left\{A\subseteq\omega:\ \forall_{F\in\cB'_n}\ F\not\subseteq A\right\}$. Put $\cB_n=\bigcup_{k\geq n}\cB'_k$ for all $n$. Then $(\cB_n)_{n\in\omega}$ is $\subseteq$-decreasing and we have:
\[
M_n=\bigcap_{k\geq n}M_k=\left\{A\subseteq\omega:\ \forall_{k\geq n}\ \forall_{F\in\cB'_k}\ F\not\subseteq A\right\}=\left\{A\subseteq\omega:\ \forall_{F\in\cB_n}\ F\not\subseteq A\right\}.
\]
\end{proof}

\begin{theorem}
\label{thm:Farah=weaklyFarah}
Let $\I$ be an ideal on $\omega$. Then $\I$ is Farah if and only if $\I$ is weakly Farah.
\end{theorem}

\begin{proof}
It is clear that every Farah ideal is weakly Farah. We will show the opposite.

Let $\I$ be a weakly Farah ideal. Then there are hereditary $\mathbf{F_\sigma}$ sets $L_n$, for $n\in\omega$, such that $\I=\bigcap_{n\in\omega}L_n$. By Proposition \ref{prop:hered_F_sigma}, there are families of non-empty finite sets $\cB^n_k$, $n,k\in\omega$, such that:
\[
L_n=\left\{A\subseteq\omega:\ \exists_{k\in\omega}\ \forall_{F\in\cB^n_k}\ F\not\subseteq A\right\},
\]
and $\cB^n_k\supseteq\cB^n_{k+1}$ for all $n,k\in\omega$. 

For every $k,n\in\omega$ and $G\in\Fin$ define:
\[
\cB^n_k(G)=\left\{F\in\cB^n_k:\ G\subsetneq F\ \land\ \min(F\setminus G)\geq k\right\}.
\]
Put also 
\[
L^G_n=\left\{A\subseteq\omega:\ \exists_{k\in\omega}\ \forall_{F\in\cB^n_k(G)}\ F\not\subseteq A\cup G\right\}.
\]
It s clear that each $L^G_n$ is $\mathbf{F_\sigma}$ and hereditary. Moreover, it is closed under finite modifications. To see this, fix $A\in L^G_n$ and $K\in\Fin$. By heredity of $L^G_n$, it is enough to show that $A\cup K\in L^G_n$. Let $k$ be such that $F\not\subseteq A\cup G$ for all $F\in\cB^n_k(G)$ and put $k'=\max\{k,1+\max K\}$. Then $\cB^n_k(G)\supseteq \cB^n_{k'}(G)$, so $F\not\subseteq (A\cup K)\cup G$ for all $F\in\cB^n_{k'}(G)$.

To finish the proof, we need to show $\I=\bigcap_{n\in\omega}\bigcap_{G\in\Fin}L^G_n$. Actually, the inclusion $\subseteq$ is easy: if $A\in\I$, $n\in\omega$ and $G\in\Fin$ are fixed, then $A\cup G\in\I$, so there is $k$ such that $F\not\subseteq A\cup G$ for all $F\in\cB^n_k$ and we get $A\in L^G_n$.

$\supseteq$: Assume that $A\in \bigcap_{n\in\omega}\bigcap_{G\in\Fin}L^G_n$. Then for every $n\in\omega$ and $G\in\Fin$ there is $k(G,n)\in\omega$ such that for every $F\in\cB^n_{k(G,n)}$ satisfying $G\subsetneq F\subseteq A$ we have $\min(F\setminus G)<k(G,n)$. We need to show that $A\in\I$. 

Inductively pick a sequence $(m_n)_{n\in\omega}$ of elements of $\omega$ such that the following conditions are satisfied for every $n\in\omega$:
\begin{itemize}
    \item[(i)] $m_0=0$ and $m_n<m_{n+1}$,
    \item[(ii)] $m_{2n+1}>\max\{k(G,n):\ G\subseteq A\cap[0,m_{2n})\}$,
    \item[(iii)] $m_{2n+2}>\max\{k(G,n):\ G\subseteq A\cap[0,m_{2n+1})\}$.
\end{itemize}
Define $B=A\cap\bigcup_{n\in\omega}[m_{2n},m_{2n+1})$ and $C=A\cap\bigcup_{n\in\omega}[m_{2n+1},m_{2n+2})$. Then $A=B\cup C$ (by (i)). Since $\I$ is an ideal, it suffices to prove that $B\in\I$ and $C\in\I$. We will only show $B\in\I$ (the second case is similar). 

Assume towards contradiction that $B\notin\I$. Then there is $n\in\omega$ such that $B\notin L_n$. Let $k_0\in\omega$ be such that $k_0>\max\{k(G,n):\ G\subseteq B\cap[0,m_{2n+1})\}$ and $\cB^n_{k_0}=\bigcup_{k\geq k_0}\cB^n_{k}$ is disjoint with $\cP([0,m_{2n+1}))$ (this is possible: by $\Fin\subseteq\I$, for every finite set $G$ there is $k\in\omega$ such that $G\notin\cB^n_k=\bigcup_{i\geq k}\cB^n_{i}$). As $B\notin L_n$, we can find $F\in B^n_{k_0}$ such that $F\subseteq B\subseteq A$. Put $G=F\cap[0,m_{2n+1})$. Then $G\subseteq B\cap[0,m_{2n+1})$, so $k(G,n)<k_0$. Hence, $F\in\cB^n_{k_0}\subseteq\cB^n_{k(G,n)}$. As $G\subsetneq F\subseteq A$, by the choice of $k(G,n)$, $\min(F\setminus G)<k(G,n)$. Since $G\subseteq A\cap[0,m_{2n+1})$, we also have $k(G,n)<m_{2n+2}$ (by (iii)). We can conclude that $\min(F\setminus G)\in[m_{2n+1},m_{2n+2})$. However, $B\cap [m_{2n+1},m_{2n+2})=\emptyset$, which contradicts $F\subseteq B$ and finishes the entire proof.
\end{proof}

%%%%%%%%%%%%%%%%%%%%%%%%%%%%%%%%%%%%%%%%%%%%%%%%%%
%%%%% SECTION
%%%%%%%%%%%%%%%%%%%%%%%%%%%%%%%%%%%%%%%%%%%%%%%%%%

\section{Counterexamples to Farah's question}

In this section we will work with submeasures and $0-1$ sequences.

Recall that a function $\phi:\mathcal{P}(\omega)\to[0,\infty]$ is called a submeasure on $\omega$ if $\phi(\emptyset)=0$, $\phi(\{n\})<\infty$ for all $n\in\N$, and 
\[
\phi(A)\leq\phi(A\cup B)\leq\phi(A)+\phi(B),
\]
for all $A,B\subseteq\omega$. A submeasure $\phi$ is lower semicontinuous (in short: lsc) if $\phi(C)=\lim_{n\to\infty}\phi(C\cap n)$ for all $C\subseteq\omega$. Given a lsc submeasure $\phi$, we can define:
\[
\Exh(\phi)=\left\{A\subseteq\omega:\ \lim_{n\to\infty}\phi(A\setminus n)=0\right\}.
\]
It is straightforward that this family is an $\mathbf{F_{\sigma\delta}}$ ideal, whenever $\phi(\omega\setminus n)\not\to 0$. By a famous result of Solecki from \cite{SoleckiExh}, an ideal is of this form if and only if it is an analytic P-ideal (for definition of a P-ideal see \cite{MR1711328}).

We will need some additional notation:
\begin{itemize}
\item If $x\in 2^\omega$ then $A_x$ stands for the set $\{x\restriction n:\ n\in\omega\}$.
\item For $A\subseteq 2^{<\omega}$ and $x\in 2^{<\omega}\cup 2^\omega$ write:
\[
L(A,x)=\left\{n\leq\lh(x):\ x\restriction n\in A\right\}
\]
(for $x\in 2^\omega$ we have $x\restriction\omega=x\notin A$, so $L(A,x)$ is always a subset of $\omega$).
\item For $x,y\in 2^{<\omega}\cup 2^\omega$ we denote:
\[
i(x,y)=\min\left(\{i\in \omega:\ x(i)\neq y(i)\}\cup\{\lh(x),\lh(y)\}\right).
\]
In particular, if $s\subseteq t\in 2^{<\omega}$ then $i(s,t)=\lh(s)$.
\end{itemize}

\begin{definition}
For each $s\in 2^{<\omega}$ and a lsc submeasure $\phi$ on $\omega$, we define the following families of subsets of $2^{<\omega}$:
\begin{equation*}
\begin{split}
X_{s,\phi}= & \left\{A\subseteq 2^{<\omega}:\ L(A,s^\frown 0_\infty)\notin\Exh(\phi)\right\}, \\
Y_{s,\phi}= & \left\{A\subseteq 2^{<\omega}:\ \exists_{k\in\omega}\ \forall_{m\in\omega}\ \exists_{s\subseteq t\in 2^{<\omega}}\ \left(\phi(L(A,t)\setminus m)\geq\frac{1}{2^k}\ \land\right.\right.\\
& \left.\left.\left(\forall_{s\subseteq t'\in 2^{\lh(t)}}\phi(L(A,t')\setminus i(t,t'))<\frac{1}{2^{k+1}}\right)\right)\right\},\\
Z_{\phi}= & \left\{A\subseteq 2^{<\omega}:\ \exists_{k\in\omega}\ \forall_{m\in\omega}\ \exists_{s\in 2^{\geq m}}\ \exists_{t_0,\ldots,t_m\in 2^{<\omega}}\ \forall_{i\leq m}\ \left(s^\frown(1)\subseteq t_i\ \land\right.\right.\\
& \left.\left.\phi\left(L(A,t_i)\setminus \max\{i(t_i,t_j):\ j\leq m,\ i\neq j\}\right)\geq\frac{1}{2^k}\right)\right\},\\
\mathcal{NF}_\phi= & \cP(2^{<\omega})\setminus \left(Z_{\phi}\cup\bigcup_{s\in 2^{<\omega}}\left(
X_{s,\phi}\cup Y_{s,\phi}\right)\right).
\end{split}
\end{equation*}
\end{definition}

\begin{lemma}
\label{lem:NF0}
Let $s\in 2^{<\omega}$ and $\phi$ be a lsc submeasure on $\omega$. Fix also $A\subseteq 2^{<\omega}$.
\begin{itemize}
    \item[(a)] If there are $k\in\omega$ and $s\subseteq x\in 2^\omega$ such that: \[
    \limsup_{m\to\infty}\phi(L(A,x)\setminus m)\geq\frac{1}{2^k}
    \]
    and for every $s\subseteq y\in 2^\omega$ we have: 
    \[
    \phi(L(A,y)\setminus i(x,y))<\frac{1}{2^{k+1}},
    \]
    then $A\in Y_{s,\phi}$.
    \item[(b)] If $A\in Y_{s,\phi}$, then there are $k\in\omega$ and $s\subseteq x\in 2^\omega$ such that: 
    \[
    \limsup_{m\to\infty}\phi(L(A,x)\setminus m)\geq\frac{1}{2^{k+1}}
    \]
    and for every $s\subseteq y\in 2^\omega$ we have: 
    \[
    \phi(L(A,y)\setminus i(x,y))<\frac{1}{2^{k+1}}.
    \]
\end{itemize}
\end{lemma}

\begin{proof}
Item (a) is straightforward. To prove (b), suppose that $A\in Y_{s,\phi}$. Then there are $k\in\omega$ and finite sequences $t_m\in 2^{<\omega}$, for $m\in\omega$, such that $s\subseteq t_m$, $\phi(L(A,t_m)\setminus m)\geq\frac{1}{2^k}$ and
\[
\phi(L(A,t')\setminus i(t_m,t'))<\frac{1}{2^{k+1}}
\]
for all $t'\in 2^{\lh(t_m)}$ with $s\subseteq t'$. The above is also satisfied for all $t'\in 2^{<\lh(t_m)}$ with $s\subseteq t'$. By passing to a subsequence, if necessary, we may assume that $\lh(t_m)<\lh(t_{m+1})$ for all $m$ (since $\phi(L(A,t_m)\setminus n)\geq\phi(L(A,t_m)\setminus m)\geq\frac{1}{2^k}$ whenever $n<m$).

If we are able to find an infinite sequence $s\subseteq t_{m_0}\subseteq t_{m_1}\subseteq\ldots$, then $x=\bigcup_{i\in\omega} t_{m_i}$ works and additionally we even get that $\limsup_{m\to\infty}\phi(L(A,x)\setminus m)\geq\frac{1}{2^{k}}$ (not only $\geq\frac{1}{2^{k+1}}$). If there is no such infinite sequence, again by passing to a subsequence, if necessary, we can assume without loss of generality that $t_n\perp t_m$ for all $n\neq m$. Observe that for each two numbers $n<m$, by the choice of $t_m$, $\phi(L(A,t_n)\setminus i(t_m,t_n))<\frac{1}{2^{k+1}}$, hence $\phi(L(A,t_n)\cap i(t_m,t_n)\setminus n)>\frac{1}{2^{k+1}}$. Define $\hat{t}_n=t_n\restriction \min\{i(t_m,t_n):\ m>n\}$ for all $n\in\omega$, so that $\phi(L(A,\hat{t}_n)\setminus n)>\frac{1}{2^{k+1}}$ and $s\subseteq\hat{t}_0\subseteq\hat{t}_1\subseteq\ldots$. Similarly as before, $x=\bigcup_{m\in\omega} \hat{t}_m$ is as required.
\end{proof}

\begin{lemma}
\label{lem:NF1}
Let $\phi$ be a lsc submeasure on $\omega$ such that $\phi(\omega\setminus n)\not\to0$.
\begin{itemize}
    \item[(a)] For every $s\in 2^{<\omega}$, the family $X_{s,\phi}$ is a co-ideal of the same complexity as $\Exh(\phi)^+$.
    \item[(b)] $Z_{\phi}$ is a $\mathbf{G_{\delta\sigma}}$ co-ideal.
    \item[(c)] $\NF_\phi$ is an $\mathbf{F_{\sigma\delta}}$ ideal.
\end{itemize}
\end{lemma}

\begin{proof}
(a): It can be easily verified that for each $s\in 2^{<\omega}$ the family $\cP(s^{<\omega})\setminus X_{s,\phi}$ is an ideal isomorphic to the ideal:
\[
\Exh(\phi)\oplus\cP(\omega)=\left\{A\subseteq 2\times\omega:\ \{n\in\omega:\ (0,n)\in A\}\in\Exh(\phi)\right\}.
\]
The ideal $\Exh(\phi)$ is either equal to $\Fin$, which has the same descriptive complexity as $\Fin\oplus\cP(\omega)$ (both ideals are $\mathbf{F_{\sigma}}$) or is isomorphic to $\Exh(\phi)\oplus\cP(\omega)$ (it is well-known and easy to show that if some infinite set belongs to an ideal $\I$, then $\I$ and $\I\oplus\cP(\omega)$ are isomorphic).

(b): It is straightforward that $Z_\phi$ is $\mathbf{G_{\delta\sigma}}$, contains only infinite sets, and that $Z_{\phi}\ni A\subseteq B$ implies $B\in Z_{\phi}$. Suppose now that $A\cup B\in Z_{\phi}$. Then there is $k\in\omega$ such that for every $m\in\omega$ we can find sequences $s_m\in 2^{\geq m}$ and $t^m_0,\ldots,t^m_m\in 2^{<\omega}$ such that $s_m\, ^\frown(1)\subseteq t^m_i$ and:
\[
\phi\left(L(A\cup B,t^m_i)\setminus \max\{i(t^m_i,t^m_j):\ j\leq m,\ i\neq j\}\right)\geq\frac{1}{2^k},
\]
for every $i\leq m$. 
Let $D$ be the set of those $m\in\omega$ such that: 
\[
\left|\left\{i\leq m:\ \phi\left(L(A,t^m_i)\setminus \max\{i(t^m_i,t^m_j):\ j\leq m,\ i\neq j\}\right)\geq\frac{1}{2^{k+1}}\right\}\right|\geq\frac{m+1}{2}.
\]
Note that if $m\in\omega\setminus D$ then the above is satisfied for the set $B$ instead of $A$.
We will now show that if $D$ is infinite then $A\in Z_\phi$ and this is witnessed by $k+1$. Indeed, for every $m\in\omega$ we can find $d\in D$ such that $\frac{d+1}{2}\geq m+1$. This means that there are pairwise distinct $i_0,\ldots,i_m\leq d$ such that:
\[
\phi\left(L(A,t^d_{i_j})\setminus \max\{i(t^d_{i_j},t^d_{i_l}):\ l\leq m,\ j\neq l\}\right)\geq\]
\[
\geq\phi\left(L(A,t^d_{i_j})\setminus \max\{i(t^d_{i_j},t^d_{i}):\ i\leq d,\ i_j\neq i\}\right)\geq\frac{1}{2^{k+1}},
\]
for all $j\leq m$. Clearly, $s_d\in 2^{\geq d}\subseteq 2^{\geq m}$ and $s_d\, ^\frown(1)\subseteq t^d_{i_j}$ for all $j\leq m$. Hence, $A\in Z_\phi$. 

Similarly we can show that if $\omega\setminus D$ is infinite then $B\in Z_\phi$.

(c): The descriptive complexity of $\NF_\phi$ is straightforward. Moreover, all finite subsets of $2^{<\omega}$ belong to $\NF_\phi$. 

Now we deal with heredity of $\NF_\phi$. We already know from previous two items that $A\in Z_{\phi}\cup \bigcup_{s\in 2^{<\omega}}X_{s,\phi}$ and $A\subseteq B$ imply $B\notin\NF_\phi$. Assume that $Y_{s,\phi}\ni A\subseteq B$ for some $s\in 2^{<\omega}$. By Lemma \ref{lem:NF0}(b), we can find $k\in\omega$ and $s\subseteq x\in 2^\omega$ such that $\limsup_{m\to\infty}\phi(L(A,x)\setminus m)\geq\frac{1}{2^{k+1}}$ and for every $s\subseteq y\in 2^\omega$ we have $\phi(L(A,y)\setminus i(x,y))<\frac{1}{2^{k+1}}$. Note that also $\lim_{m\to\infty}\phi(L(B,x)\setminus m)\geq\frac{1}{2^{k+1}}$.

Let $Y$ be the set of all $y\in 2^\omega$ satisfying:
    \[
    \phi(L(B,y)\setminus i(x,y))\geq\frac{1}{2^{k+2}}.
    \]
There are three possible cases:
\begin{itemize}
    \item If $x\in\mathbb{Q}(2^\omega)$ then $B\in X_{x\restriction n,\phi}$, so $B\notin\NF_\phi$.
    \item If we can find $n>\lh(s)$ such that $[x\restriction n]\cap Y=\emptyset$, then Lemma \ref{lem:NF0}(a) gives us $B\in Y_{x\restriction n,\phi}$, so again $B\notin\NF_\phi$.
    \item If the above two cases do not hold, then $x\notin\mathbb{Q}(2^\omega)$ and for every $n\in\omega$ we are able to find $y_n\in Y$ such that $x\restriction n=y_n\restriction n$. By passing to a subsequence, without loss of generality we may assume that $(i(x,y_n))_{n\in\omega}$ is strictly increasing. By lower semicontinuity of $\phi$, for each $n$ we can find a finite sequence $y_n\restriction n\subseteq t_n\subseteq y_n$ such that $\phi(L(B,t_n)\setminus i(x,t_n))\geq\frac{1}{2^{k+3}}$. We will show that $B\in Z_\phi$. Indeed, fix $m\in\omega$ and find $l>m$ such that $x(l)=1$ (which exists by $x\notin\mathbb{Q}(2^\omega)$). Then for $s'=x\restriction l\in 2^{\geq m}$ and $t'_0=t_{l+1},t'_1=t_{l+2},\ldots,t'_m=t_{l+1+m}$ we have $s'^\frown (1)=x\restriction (l+1)\subseteq t'_i$ and:
    \[
    \phi\left(L(B,t'_i)\setminus \max\{i(t'_i,t'_j):\ j\leq m,\ i\neq j\}\right)\geq\phi\left(L(B,t'_i)\setminus i(x,t'_i)\right)\geq\frac{1}{2^{k+3}}
    \]
    for all $i\leq m$. Hence, $B\in Z_\phi$. 
\end{itemize}
This shows heredity of $\NF_\phi$.

Lastly, we show that $A\cup B\notin\NF_\phi$ implies $A\notin\NF_\phi$ or $B\notin\NF_\phi$. By items (a) and (b), this is clear whenever $A\cup B\in Z_{\phi}\cup \bigcup_{s\in 2^{<\omega}}X_{s,\phi}$. Therefore, suppose that $A\cup B\in Y_{s,\phi}$ for some $s$. By Lemma \ref{lem:NF0}(b), we can find $k\in\omega$ and $s\subseteq x\in 2^\omega$ such that $\limsup_{m\to\infty}\phi(L(A\cup B,x)\setminus m)\geq\frac{1}{2^{k+1}}$ and for every $s\subseteq y\in 2^\omega$ we have: 
    \[
    \phi(L(A\cup B,y)\setminus i(x,y))<\frac{1}{2^{k+1}}.
    \]
The properties of submeasures give us that either $\limsup_{m\to\infty}\phi(L(A,x)\setminus m)\geq\frac{1}{2^{k+2}}$ or $\limsup_{m\to\infty}\phi(L(B,x)\setminus m)\geq\frac{1}{2^{k+2}}$. Without loss of generality we will assume that the former case holds. Then $A\cap A_x\in Y_{s,\phi}$ (by Lemma \ref{lem:NF0}(a)) and using heredity of $\NF_\phi$, which we have already proved above, $A\notin\NF_{\phi}$. 
\end{proof}

\begin{proposition}
\label{prop:NF}
Let $\phi$ be a lsc submeasure on $\omega$ such that $\phi(\omega\setminus n)\not\to0$.
Then $\NF_\phi$ and $\Exh(\phi)$ are $\leq_K$-equivalent. In particular, $\NF_\phi$ is tall if and only if $\Exh(\phi)$ is tall.
\end{proposition}

\begin{proof}
The function $f:\omega\to 2^{<\omega}$ given by $f(n)=0_n$ witnesses $\NF_\phi\leq_K\Exh(\phi)$. Indeed, it follows directly from the definition of $X_{\emptyset,\phi}$ that if $A\in\Exh(\phi)^+$ then $f[A]\in X_{\emptyset,\phi}\subseteq\NF_\phi^+$.

On the other hand, $g:2^{<\omega}\to\omega$ given by $g(s)=\lh(s)$ witnesses $\Exh(\phi)\leq_K\NF_\phi$, since $A\in\NF_\phi^+$ implies $g[A]\in\Exh(\phi)^+$.

It is a well-known and easy to see fact that any ideal $\I$ is not tall if and only if $\I\leq_K\Fin$. Hence, the "In particular" part follows immediately from the first part.
\end{proof}

\begin{lemma}
\label{lem:NF2}
For every lsc submeasure $\phi$ on $\omega$ such that $\phi(\omega\setminus n)\not\to0$, the ideal $\NF_\phi$ is not Farah.
\end{lemma}

\begin{proof}
Actually, we will show that $\NF_\phi$ is not weakly Farah (which by Theorem \ref{thm:Farah=weaklyFarah} is equivalent to not being Farah). Assume towards contradiction that there are hereditary $\mathbf{F_\sigma}$ sets $L_n$ such that $\NF_\phi=\bigcap_{n\in\omega}L_n$. By Proposition \ref{prop:hered_F_sigma}, there are families of non-empty finite subsets $\cB^n_k$ of $2^{<\omega}$, $n,k\in\omega$, such that:
\[
L_n=\left\{A\subseteq 2^{<\omega}:\ \exists_{k\in\omega}\ \forall_{F\in\cB^n_k}\ F\not\subseteq A\right\}.
\]

Since $\phi(\omega\setminus n)\not\to0$, there are $\delta>0$ and a partition of $\omega$ into finite intervals $I_n$, $n\in\omega$, such that $\phi(I_n)>\delta$ and $1+\max I_n=\min I_{n+1}$ for all $n$.

For $x\in 2^\omega$ and $B\in[\omega]^\omega$ denote $C(x,B)=\bigcup_{n\in B}\{x\restriction i:\ i\in I_n\}$. By Lemma \ref{lem:NF0}(a), each set of the form $C(x,B)$ belongs to $Y_{\emptyset,\phi}$, so also to $\NF^+_\phi$. Thus, for every $x\in 2^\omega$ and $B\in[\omega]^\omega$ there is $n(x,B)$ such that $C(x,B)\notin L_{n(x,B)}$. Since $[\omega]^\omega$ is a dense $\mathbf{G_\delta}$ subset of $\cP(\omega)$, by the Baire category theorem, for each $x\in 2^\omega$ there are $n_x\in\omega$ and a basic set $U_x\subseteq\mathcal{P}(\omega)$ such that the set of those $B\in[\omega]^\omega$, for which $n(x,B)=n_x$ is dense in $U_x$. Without loss of generality we can assume that each $U_x$ is of the form $\{A\subseteq\omega:\ A\cap p_x=F_x\}$ for some finite set $F_x\subseteq\omega$ and $p_x\in\omega$. Now, using the Baire category theorem once again, we can find $n,p\in\omega$, $F\in\Fin$ and a finite sequence $s\in 2^{<\omega}$ such that $\{x\in 2^\omega:\ n=n_x\ \land\ F=F_x\ \land\ p=p_x\}$ is dense in $[s]$. Without loss of generality we will assume that $\text{lh}(s)\geq \max I_{p}$.

For each $k\in\omega$ find:
\begin{itemize}
    \item $x_k\in [s]$ such that $n_{x_k}=n$, $F_{x_k}=F$, $p_{x_k}=p$ and $s^\frown 0_k^\frown(1)\subseteq x_k$,
    \item $B_k\in U_{x_k}$ such that $n(x_k,B_k)=n_{x_k}$ and:
    \[
    C(x_k,B_k)\cap \{x\restriction i:\ i\leq \lh(s^\frown 0_k^\frown(1))\}=\bigcup_{n\in F}\{x_k\restriction i:\ i\in I_n\},
    \]
    \item $F_k\in\mathcal{B}^n_k$ such that $F_k\subseteq C(x_k,B_k)$ (which exists by $C(x_k,B_k)\notin L_n$).
\end{itemize}
Define $A=\bigcup_{k\in\omega}F_k$. Then $A\notin L_n$, so $A\in\NF_\phi^+$. On the other hand, we will show that $A\in\NF_\phi$. Indeed, $A\notin\bigcup_{t\in 2^{<\omega}}(X_{t,\phi}\cup Y_{t,\phi})$ as $A\cap A_x$ is finite for all $x\in 2^\omega$ (cf. Lemma \ref{lem:NF0}(b)). Moreover, note that each $A\cap\{t'\in 2^{<\omega}:\ t^\frown(1)\subseteq t'\}$, with $t\in 2^{\geq \lh(s)}$, is either empty or can be covered by a single set of the form $A_x$. Hence, $A\notin Z_\phi$. 
\end{proof}

\begin{theorem}\
\label{thm:NF}
\begin{itemize}
    \item[(a)] There is an $\mathbf{F_{\sigma\delta}}$ non-tall ideal that is not Farah.
    \item[(b)] There is a $\leq_K$-antichain of length $\cc$ consisting of $\mathbf{F_{\sigma\delta}}$ tall ideals that are not Farah.
\end{itemize}
\end{theorem}

\begin{proof}\
(a): Let $\phi$ be any lsc submeasure on $\omega$ such that $\Exh(\phi)=\Fin$ (e.g. $\phi(A)=|A|$ for all $A\subseteq\omega$). It follows from Lemmas \ref{lem:NF1}(c) and \ref{lem:NF2} that $\NF_\phi$ is an $\mathbf{F_{\sigma\delta}}$ non-Farah ideal. Moreover, $\NF_\phi$ is not tall by Proposition \ref{prop:NF}.

(b): By \cite[Corollary 1]{MezaGuzman}, there is a family $\cA$ of cardinality $\cc$ consisting of pairwise $\leq_K$-incomparable analytic P-ideals on $\omega$, that is, ideals of the form $\Exh(\phi)$ (alternatively one can use \cite[Theorem 3]{b} together with \cite[Proposition 1.1]{simpleden}). Since all non-tall ideals are $\leq_K$-equivalent, without loss of generality we can assume that all members of $\cA$ are tall. For each $\I\in\cA$ let $\phi_\I$ be a lsc submeasure on $\omega$ such that $\I=\Exh(\phi_\I)$. Applying Lemmas \ref{lem:NF1}(c) and \ref{lem:NF2} to each $\I\in\cA$, we get that $\NF_{\phi_\I}$ is an $\mathbf{F_{\sigma\delta}}$ non-Farah ideal. Moreover, by Proposition \ref{prop:NF}, if $\I,\J\in\cA$ are distinct, then $\NF_{\phi_\I}$ and $\NF_{\phi_\J}$ are $\leq_K$-incomparable and tall.
\end{proof}

%%%%%%%%%%%%%%%%%%%%%%%%%%%%%%%%%%%%%%%%%%%%%%%%%%
%%%%% SECTION
%%%%%%%%%%%%%%%%%%%%%%%%%%%%%%%%%%%%%%%%%%%%%%%%%%

\section{Sets of ideal limit points}

\begin{definition}
\label{def:gra}
Let $\I$ be an ideal on $\omega$ and $\cA=\{A_s:\ s\in\omega^{<\omega}\}$ be a system of sets such that for all $s\in\omega^{<\omega}$ and $i,j\in\omega$ the following conditions hold:
\begin{itemize}
    \item $A_{s^\frown(i)}\subseteq A_s$,
    \item $A_{s^\frown(i)}\cap A_{s^\frown(j)}=\emptyset$ whenever $i\neq j$.
\end{itemize}
The game $G(\I,\cA)$ is played by two players in the following way: in the initial move Player I plays $s_0\in\omega^{<\omega}$ and Player II responses with $s_0\subseteq t_0\in\omega^{<\omega}$ and a finite set $F_0\subseteq A_{s_0}$. In the $n$th move, $n>0$, Player I plays $s_n\in\omega^{<\omega}$ such that $t_{n-1}\subsetneq s_n$ and then Player II responses with $s_n\subseteq t_n\in\omega^{<\omega}$ and a finite set $F_n\subseteq A_{s_n}$. Player I is declared the winner if $\bigcup_{n\in\omega}F_n\in\I$. Otherwise Player II is the winner.
\end{definition}

\begin{lemma}
\label{lem:gra}
Let $\I$ be an ideal on $\omega$ and $\cA$ be as in Definition \ref{def:gra}. Assume that for every $x\in\omega^\omega$ there is $A\in\I^+$ such that $A\setminus A_{x\restriction n}$ is finite for all $n$. Then Player I does not have a winning strategy in $G(\I,\cA)$.
\end{lemma}

\begin{proof}
Assume towards contradiction that Player I has a winning strategy $\sigma$. Let $f:([\omega]^{<\omega})^{<\omega}\to\omega$ be such a bijection that $s\subseteq t$ implies $f(s)\leq f(t)$ for all $s,t,\in ([\omega]^{<\omega})^{<\omega}$ (note that this condition implies $f(\emptyset)=0$).

Recursively define a sequence $(u_n)_{n\in\omega}$ of elements of $\omega^{<\omega}$ by putting $u_0=u_{f(\emptyset)}=\sigma(\emptyset)$ and:
\begin{equation*}
\begin{split}
    u_{n}= & \sigma\left(F_0\cap A_{u_{f(\emptyset)}},u_{f(F_0)-1},F_1\cap A_{u_{f(F_0)}},u_{f(F_0,F_1)-1},\ldots,\right.\\
    & \left.F_k\cap A_{u_{f(F_0,\ldots,F_{k-1})}},u_{f(F_0,F_1,\ldots,F_{k})-1}\right),
\end{split}
\end{equation*}
for all $n$, where $(F_0,\ldots,F_k)=f^{-1}(n)$. We need to inductively show that all $u_n$'s are well-defined and that $u_n\subsetneq u_{n+1}$ for all $n$. Assume that all $u_i$, for $i<n=f(F_0,\ldots,F_k)$, are well-defined and satisfy $u_0\subsetneq u_1\subsetneq\ldots\subsetneq u_{n-1}$. Observe that:
\begin{itemize}
    \item $F_k\cap A_{u_{f(F_0,\ldots,F_{k-1})}}\subseteq A_{u_{f(F_0,\ldots,F_{k-1})}}$,
    \item $f(F_0,F_1,\ldots,F_{k-1})\leq f(F_0,F_1,\ldots,F_{k})-1<f(F_0,F_1,\ldots,F_{k})=n$ and, consequently, $u_{f(F_0,F_1,\ldots,F_{k})-1}\supseteq u_{f(F_0,F_1,\ldots,F_{k-1})}$ (by the inductive assumption).
\end{itemize}
This shows that $F_k\cap A_{u_{f(F_0,\ldots,F_{k-1})}}$ and $u_{f(F_0,F_1,\ldots,F_{k})-1}$ are legal plays of Player II after Player I played:
\begin{equation*}
\begin{split}
    u_{f(F_0,\ldots,F_{k-1})}= & \sigma\left(F_0\cap A_{u_{f(\emptyset)}},u_{f(F_0)-1},F_1\cap A_{u_{f(F_0)}},u_{f(F_0,F_1)-1},\ldots,\right.\\
    & \left.F_{k-1}\cap A_{u_{f(F_0,\ldots,F_{k-2})}},u_{f(F_0,F_1,\ldots,F_{k-1})-1}\right).
\end{split}
\end{equation*}
Hence, $u_n$ is well-defined. Moreover, it follows from the rules of our game that $u_n\supsetneq u_{f(F_0,F_1,\ldots,F_{k})-1}=u_{n-1}$.

Let $x=\bigcup_{n\in\omega}u_n\in\omega^\omega$ (which is well-defined thanks to the condition $u_n\subsetneq u_{n+1}$). By the assumptions of the whole Lemma, there is $A\in\I^+$ such that $A\setminus A_{x\restriction n}$ is finite for all $n$. Then also $A\setminus A_{u_n}$ is finite for all $n$. Without loss of generality we may assume that $A\subseteq A_{u_0}$. Let $(P_n)_{n\in\omega}$ be any partition of $\bigcap_{n\in\omega}A_{x\restriction n}$ into finite sets (if $\bigcap_{n\in\omega}A_{x\restriction n}$ is empty, then put $P_n=\emptyset$ for all $n$).

Now we recursively define two sequences of finite sets $(G_n)_{n\in\omega}$ and $(H_n)_{n\in\omega}$ by:
\begin{itemize}
    \item $G_0=(A\setminus A_{u_1})\cup P_0$,
    \item $H_0=(A\cap A_{u_1})\setminus A_{u_{f(G_0)}}$,
    \item $G_{n+1}=((A\cap A_{u_{f(G_0,\ldots,G_n)}})\setminus A_{u_{f(H_0,\ldots,H_n)}})\cup P_{n+1}$ for all $n$,
    \item $H_{n+1}=(A\cap A_{u_{f(H_0,\ldots,H_n)}})\setminus A_{u_{f(G_0,\ldots,G_n,G_{n+1})}}$ for all $n$.
\end{itemize}

Thanks to the choice of the bijection $f$, the sequences $(f(G_0,\ldots,G_n))_{n\in\omega}$ and $(f(H_0,\ldots,H_n))_{n\in\omega}\in\omega^\omega$ are increasing. Thus, we have:
\[
\I^+\ni A=\bigcup_{n\in\omega}G_n\cup\bigcup_{n\in\omega}H_n.
\]
Hence, either $\bigcup_{n\in\omega}G_n\in\I^+$ or $\bigcup_{n\in\omega}H_n\in\I^+$. We will obtain a contradiction by showing that both $\bigcup_{n\in\omega}G_n$ and $\bigcup_{n\in\omega}H_n$ belong to $\I$. 

Consider the play of the game $G(\I,\cA)$, in which:
\begin{itemize}
    \item Player I starts with $\sigma(\emptyset)=u_0$,
    \item Player II in his first move plays $G_0\subseteq A_{u_0}$ and $u_{f(G_0)-1}\supseteq u_0$,
    \item Player I in his second move plays $\sigma(G_0,u_{f(G_0)-1})=u_{f(G_0)}$,
    \item Player II in his second move plays $G_1\subseteq A_{u_{f(G_0)}}$ and $u_{f(G_0,G_1)-1}\supseteq u_{f(G_0)}$,
    \item Player I in his third move plays $\sigma(G_0,u_{f(G_0)-1},G_1,u_{f(G_0,G_1)-1})=u_{f(G_0,G_1)}$,
    \item Player II in his third move responses with $G_2\subseteq A_{u_{f(G_0,G_1)}}$ together with $u_{f(G_0,G_1,G_2)-1}\supseteq u_{f(G_0,G_1)}$,
    \item etc.
\end{itemize}
This is a legal play, in which Player I follows his winning strategy $\sigma$, so $\bigcup_{n\in\omega}G_n\in\I$. Similarly one can show that $\bigcup_{n\in\omega}H_n\in\I$.
\end{proof}

\begin{theorem}
If $\I$ is an $\bf{G_{\delta\sigma\delta}}$ ideal on $\omega$ and $X$ is an uncountable Polish space, then $\mathcal{L}_X(\I)$ is equal either to the family of all closed subsets of $X$ or to the family of all $\bf{F_{\sigma}}$ subsets of $X$.
\end{theorem}

\begin{proof}
By \cite[Corollary 10.6]{FKL24} it suffices to show that there is no family $\cA=\{A_s:\ s\in 2^{<\omega}\}$ of subsets of $\omega$ satisfying for all $s\in 2^{<\omega}$:
\begin{itemize}
    \item[(i)] $A_s\in\I^+$,
    \item[(ii)] $A_{s^\frown(0)}\cap A_{s^\frown(1)}=\emptyset$,
    \item[(iii)] $A_{s^\frown(0)}\cup A_{s^\frown(1)}\subseteq A_s$,
    \item[(iv)] $\cB_\I(\cA)=\mathbb{Q}(2^\omega)$, where:
    \[
    \cB_\I(\cA)=\left\{x\in 2^\omega:\ \forall_{B\in\I^+}\ \exists_{n\in\omega}\ B\setminus A_{x\restriction n}\in\Fin^+\right\}.
    \]
\end{itemize}
Thus, suppose towards contradiction that there is a family $\cA=\{A_s:\ s\in 2^{<\omega}\}$ as above. 

Define $f:\omega^{<\omega}\to 2^{<\omega}$ by: 
\[
f(s)=0_{s(1)}\, ^{\frown}(1)^{\frown}0_{s(2)}\, ^\frown (1)^\frown\ldots^\frown0_{s(k)}\, ^\frown (1),
\]
for all $s=(s(0),s(1),\ldots,s(k))\in\omega^{<\omega}$. In particular, $f(s^\frown(i))=f(s)^\frown 0_i^\frown(1)$. Observe that for each $s\in\omega^{<\omega}$ and distinct $i,j\in\omega$:
\begin{itemize}
    \item $A_{f(s^\frown(i))}\subseteq A_{f(s)}$,
    \item $A_{f(s^\frown(i))}\cap A_{f(s^\frown(j))}=\emptyset$.
\end{itemize}
Moreover, since $\cB_\I(\cA)=\mathbb{Q}(2^\omega)$, for every $x\in\omega^\omega$ there is $A\in\I^+$ such that $A\setminus A_{f(x\restriction n)}$ is finite for all $n$. Hence, it makes sense to consider the game $G(\I,\{A_{f(s)}:\ s\in\omega^{<\omega}\})$ and, by Lemma \ref{lem:gra}, Player I cannot have a winning strategy in this game. However, since $\I$ is Borel, this game is determined (by Martin's Theorem on Borel Determinacy). Thus, there is a winning strategy for Player II. Let $\sigma_0(s_0,\ldots,s_k)\in\Fin$ and $\sigma_1(s_0,\ldots,s_k)\in \omega^{<\omega}$ denote the moves of Player II according to this strategy after Player I played in his first $k+1$ moves $s_0\subsetneq \ldots\subsetneq s_k$ (and Player II followed his strategy in the previous moves).

Let $\{U_{k,n,m}:\ k,n,m\in\omega\}$ be a family of open sets such that:
\[
\I=\bigcap_{k\in\omega}\bigcup_{n\in\omega}\bigcap_{m\in\omega}U_{k,n,m}
\]
and $U_{k,n,m+1}\subseteq U_{k,n,m}$ for all $k,n,m\in\omega$. Fix also any bijection $g:\omega\to\Fin$ such that $g(0)=\emptyset$. 

We need to perform a recursive construction of a sequence $(n_k)_{k\in\omega}$ of natural numbers together with a decreasing sequence:
\[
B^0_0\supseteq B^1_0\supseteq B^1_1\supseteq\ldots\supseteq B^n_0\supseteq B^n_1\supseteq\ldots\supseteq B^n_n\supseteq B^{n+1}_0\supseteq\ldots
\]
of non-empty basic subsets of $\cP(2^{<\omega})$ such that $B^m_k\subseteq U_{k,n_k,m}$ for all $k\in\omega$ and $m\geq k$. Additionally, we will build a sequence $(l_k)_{k\in\omega}$ defining a play of the game $G(\I,\{A_{f(s)}:\ s\in\omega^{<\omega}\})$, in which Player II follows his winning strategy, while Player I starts with $(l_0)\in\omega^{<\omega}$ and after that in each of his moves extends the previous move $t_k\in\omega^{<\omega}$ of Player II by playing $t_k\, ^\frown(l_{k+1})\in\omega^{<\omega}$ (in particular, $t_k=\sigma_1((l_0),t_0\, ^\frown(l_1),\ldots,t_{k-1}\, ^\frown(l_k))$). The basic sets $B^m_k$ will be connected to the sequence $(l_k)_{k\in\omega}$ -- in the construction we will require that $G^k_0\supseteq \sigma_0((l_0),t_0\, ^\frown(l_1),\ldots,t_{k-1}\, ^\frown(l_k))$, where $G^m_k$ denotes the finite set such that $B^m_k=\{A\subseteq 2^{<\omega}:\ G^m_k\subseteq A,\ H^m_k\cap A=\emptyset\}$.

Before performing the construction, we show how the rest of the proof follows. Assuming that the above described objects are already defined, observe that the finite sets $G^m_k$ and $H^m_k$ must satisfy:
\[
G^0_0\subseteq G^1_0\subseteq G^1_1\subseteq\ldots\subseteq G^n_0\subseteq G^n_1\subseteq\ldots\subseteq G^n_n\subseteq G^{n+1}_0\subseteq\ldots,
\]
\[
H^0_0\subseteq H^1_0\subseteq H^1_1\subseteq\ldots\subseteq H^n_0\subseteq H^n_1\subseteq\ldots\subseteq H^n_n\subseteq H^{n+1}_0\subseteq\ldots.
\]
Define $A=\bigcup_{n\in\omega} G^n_n$. Since $t_k=\sigma_1((l_0),t_0\, ^\frown(l_1),\ldots,t_{k-1}\, ^\frown(l_k))$ and:
\[
G^k_0\supseteq F_k=\sigma_0((l_0),t_0\, ^\frown(l_1),\ldots,t_{k-1}\, ^\frown(l_k)),
\]
the set $A$ contains the result $\bigcup_{k\in\omega}F_k$ of a play of $G(\I,\{A_{f(s)}:\ s\in\omega^{<\omega}\})$, in which Player II followed his winning strategy. Hence, $\I^+\ni \bigcup_{k\in\omega}F_k\subseteq A$. On the other hand, 
\[
A\in\bigcap_{k\in\omega}\bigcup_{n\in\omega}\bigcap_{m\in\omega}U_{k,n,m}=\I,
\]
as for every $k\in\omega$ there is $n_k$ such that for every $m\geq k$ we have $A\in B^m_k\subseteq U_{k,n_k,m}$ ($A\notin B^m_k$ would mean that for some $i$ we have $\emptyset\neq G^i_i\cap H^m_k\subseteq G^{j}_{j}\cap H^{j}_{j}$, where $j=\max\{k,m,i\}$, which is contrary to $B^j_j\neq\emptyset$) and $(U_{k,n_k,m})_{m\in\omega}$ is decreasing. This contradiction will finish the whole proof.

\medskip
{\bf{First step of the construction:}} 

Since $0_\infty\in B_\I(\cA)$, for every $x\in\omega^\omega$ the set:
\[
E_0(x)=\bigcup_{i\in\omega}(g(x(i))\cap A_{f((i))})
\]
has to belong to $\I$ (recall that $A_{f((i))}=A_{0_i^\frown(1)}$). Thus, for each $x\in\omega^\omega$ we can find $n_0(x)\in\omega$ such that $E_0(x)\in\bigcap_{m\in\omega}U_{0,n_0(x),m}$. Using Baire category theorem, find $w^0_0\in\omega^{<\omega}$ and $n_0\in\omega$ such that $\{x\in\omega^\omega:\ n_0(x)=n_0\}$ is dense in $[w^0_0]$. Let $l_0=\lh(w^0_0)$ and put $F_0=\sigma_0((l_0))\subseteq A_{f((l_0))}=A_{0_{l_0}^\frown(1)}$ and $t_0=\sigma_1((l_0))\in\omega^{<\omega}$ (so that $(l_0)\subseteq t_0$). Note that $x\in[w^0_0]$ if and only if:
\[
E_0(x)\cap \bigcup_{i<l_0}A_{f((i))}=\bigcup_{i<l_0}(g(w^0_0(i))\cap A_{f((i))}).
\]

Find any $x^0_0\in [w^0_0]$ such that $n_0(x^0_0)=n_0$ and $g(x^0_0(l_0))=F_0$. Then we can find a basic set $B^0_0=\{A\subseteq 2^{<\omega}:\ G^0_0\subseteq A,\ H^0_0\cap A=\emptyset\}$, for some finite sets $G^0_0$ and $H^0_0$, such that $E_0(x^0_0)\in B^0_0\subseteq U_{0,n_0,0}$. Observe that:
\[
E_0(x^0_0)\in B^0_0\cap\left\{A\subseteq 2^{<\omega}:\ F_0\cup \bigcup_{i<l_0}(g(w^0_0(i))\cap A_{f((i))})\subseteq A\right\}\subseteq B^0_0\subseteq U_{0,n_0,0}.
\]
Hence, without loss of generality we may assume that:
\[
G^0_0\supseteq F_0\cup \bigcup_{i<l_0}(g(w^0_0(i))\cap A_{f((i))}).
\]
Then we have: 
\[
\left( F_0\cup \bigcup_{i<l_0}(g(w^0_0(i))\cap A_{f((i))})\right)\cap H^0_0\subseteq G^0_0\cap H^0_0=\emptyset
\]
(as $B^0_0\neq\emptyset$) and:
\[
G^0_0\setminus\left( F_0\cup \bigcup_{i<l_0}(g(w^0_0(i))\cap A_{f((i))})\right)\subseteq\bigcup_{i>l_0} A_{f((i))},
\]
(by $x^0_0\in [w^0_0]$, $g(x^0_0(l_0))=F_0$ and $E_0(x^0_0)\in B^0_0$). This finishes the initial step of our construction. 

\medskip
{\bf{Second step of the construction:}} 

Since $0_\infty\in \cB_\I(\cA)$ and $f(t_0)\, ^\frown 0_\infty\in\cB_\I(\cA)$, for every $\bar{x}=(x_0,x_1)\in\omega^\omega\times\omega^\omega$ the set:
\[
E_1(\bar{x})=F_0\cup \bigcup_{i\in\omega\setminus\{l_0\}}(g(x_0(i))\cap A_{f((i))})\cup \bigcup_{i\in\omega}(g(x_1(i))\cap A_{f(t_0^\frown(i))})
\]
has to belong to $\I$. For each $\bar{x}\in\omega^\omega\times\omega^\omega$ find $n_1(\bar{x})\in\omega$ such that $E_1(\bar{x})\in\bigcap_{m\in\omega}U_{1,n_1(\bar{x}),m}$. Note that $E_1(\bar{x})$ does not depend on $x_0(l_0)$, so we can assume that $n_1(\bar{x})$ also does not depend on $x_0(l_0)$, that is, $n_1(x_0,x_1)=n_1(x'_0,x'_1)$ whenever $x_1=x'_1$ and $x_0\restriction \omega\setminus\{l_0\}=x'_0\restriction \omega\setminus\{l_0\}$.

By the observations ending first step of the construction, it is possible to find $u^{1}_0,u^{1}_1\in\omega^{<\omega}$ such that $w^0_0\subseteq u^{1}_0$ and $\bar{x}\in[u^{1}_0]_{l_0}\times[u^{1}_1]$ implies $E_1(\bar{x})\in B^0_0$, where for $w\in\omega^{<\omega}$ and $l\in\omega$ we denote:
\[
[w]_{l}=\{x\in\omega^\omega:\ x(i)=w(i)\text{ for all }i<\lh(w), i\neq l\}.
\]
Using Baire category theorem, find $n_1\in\omega$ such that $\{\bar{x}\in\omega^\omega\times\omega^\omega:\ n_1(\bar{x})=n_1\}$ is somewhere dense in $[u^{1}_0]_{l_0}\times[u^{1}_1]$. Since $n_1(\bar{x})$ does not depend on $x_0(l_0)$, the above means that there are $w^1_0,w^1_1\in\omega^{<\omega}$ such that $[w^1_0]_{l_0}\times[w^1_1]\subseteq[u^{1}_0]_{l_0}\times[u^{1}_1]$ and $\{\bar{x}\in\omega^\omega\times\omega^\omega:\ n_1(\bar{x})=n_1\}$ is dense in $[w^1_0]_{l_0}\times[w^1_1]$. Let $l_1=\lh(w^1_1)$ and put $t_1=\sigma_1((l_0),t_0\, ^\frown(l_1))\in\omega^{<\omega}$ (so that $(l_0)\subseteq t_0\subseteq t_0\, ^\frown(l_1)\subseteq t_1$) and $F_1=\sigma_0((l_0),t_0\, ^\frown(l_1))$ (in particular, $F_1\subseteq A_{f(t_0\, ^\frown(l_1))}=A_{f(t_0)^\frown 0_{l_1}\,^\frown (1)}\subseteq A_{f((l_0))}$). 

Find $x^1_0\in[w^1_0]_{l_0}\subseteq[w^0_0]$ such that $n_0(x^1_0)=n_0$ and:
\[
g(x^1_0(l_0))=F_0\cup F_1\cup \bigcup_{i<l_1}(g(w^1_1(i))\cap A_{f(t_0^\frown (i))}).
\]
Observe that $E_0(x^1_0)=E_1(x^1_0,w^1_1 \,^\frown(g^{-1}(F_1))^\frown 0_\infty)\in B^0_0$ (since $g(0)=\emptyset$ and $(x^1_0,w^1_1 \,^\frown(g^{-1}(F_1))^\frown 0_\infty)\in [u^{1}_0]_{l_0}\times[u^{1}_1]$). Then we can find a basic set $B^1_0=\{A\subseteq 2^{<\omega}:\ G^1_0\subseteq A,\ H^1_0\cap A=\emptyset\}$, for finite $G^1_0$ and $H^1_0$, such that $E_0(x^1_0)\in B^1_0\subseteq U_{0,n_0,1}\cap B^0_0$. Similarly as in the first step, without loss of generality we may assume that:
\[
G^1_0\supseteq g(x^1_0(l_0))\cup \bigcup_{i<\lh(w^1_0),i\neq l_0}(g(w^1_0(i))\cap A_{f((i))}) \supseteq F_1\cup G^0_0.
\]
Find also any $\bar{x}^1_1=(x^1_{1,0},x^1_{1,1})\in [w^{1}_0]_{l_0}\times[w^{1}_1]$ such that $n_1(\bar{x}^1_1)=n_1$, $g(x^1_{1,1}(l_1))=F_1$ and $E_1(\bar{x}^1_1)\in B^1_0$ (note that this is possible). Then there is $B^1_1=\{A\subseteq 2^{<\omega}:\ G^1_1\subseteq A,\ H^1_1\cap A=\emptyset\}$, $G^1_1$ and $H^1_1$ finite, such that $E_1(\bar{x}^1_1)\in B^1_1\subseteq U_{1,n_1,1}\cap B^1_0$. 

We have: 
\[
\left( g(x^1_0(l_0))\cup \bigcup_{i<\lh(w^1_0),i\neq l_0}(g(w^1_0(i))\cap A_{f((i))})\right)\cap H^1_1\subseteq G^1_0\cap H^1_1\subseteq G^1_1\cap H^1_1=\emptyset,
\]
\[
G^1_1\setminus\left( g(x^1_0(l_0))\cup \bigcup_{i<\lh(w^1_0),i\neq l_0}(g(w^1_0(i))\cap A_{f((i))})\right)\subseteq\bigcup_{i\geq\lh(w^1_0)} A_{f((i))}\cup \bigcup_{i>l_1} A_{f(t_0^\frown(i))}.
\]

\medskip
{\bf{Third step of the construction:}} 

Now, since $0_\infty,f(t_0)\, ^\frown 0_\infty,f(t_1)\, ^\frown 0_\infty\in \cB_\I(\cA)$, for every $\bar{x}=(x_0,x_1,x_2)\in(\omega^\omega)^3$ the set:
\[
E_2(\bar{x})=   F_0\cup F_1\cup \bigcup_{i\in\omega\setminus\{l_0\}}(g(x_0(i))\cap A_{f((i))})\cup 
\]
\[
\bigcup_{i\in\omega\setminus\{l_1\}}(g(x_1(i))\cap A_{f(t_0^\frown (i))})\cup \bigcup_{i\in\omega}(g(x_2(i))\cap A_{f(t_1^\frown (i))})
\]
has to belong to $\I$, so there is $n_2(\bar{x})\in\omega$ such that $E_2(\bar{x})\in\bigcap_{m\in\omega}U_{2,n_2(\bar{x}),m}$. Since $E_2(\bar{x})$ does not depend on $x_0(l_0)$ and $x_1(l_1)$, we can assume that $n_2(\bar{x})$ also does not depend on those values.

Find $u^{2}_0,u^{2}_1,u^{2}_2\in\omega^{<\omega}$ such that $w^{1}_0\subseteq u^{2}_0$, $w^{1}_1\subseteq u^{2}_1$ and $\bar{x}\in [u^{2}_0]_{l_0}\times[u^{2}_1]_{l_1}\times[u^{2}_2]$ implies $E_2(\bar{x})\in B^1_1$ (which is possible by the observations at the end of the previous step). By Baire category theorem, there are $w^2_0,w^2_1,w^2_2\in\omega^{<\omega}$ and $n_2\in\omega$ such that $\{\bar{x}\in(\omega^\omega)^3:\ n_2(\bar{x})=n_2\}$ is dense in $[w^2_0]_{l_0}\times[w^2_1]_{l_1}\times[w^2_2]\subseteq [u^{2}_0]_{l_0}\times[u^{2}_1]_{l_1}\times[u^{2}_2]$. Let $l_2=\text{lh}(w^2_2)$ and put $F_2=\sigma_0((l_0),t_0\, ^\frown(l_1),t_1\, ^\frown(l_2))\subseteq A_{f(t_1^\frown(l_2))}$ and $t_2=\sigma_1((l_0),t_0\, ^\frown(l_1),t_1\, ^\frown(l_2))$.

Let $x^2_0\in[w^2_0]_{l_0}\subseteq[w^0_0]$ be such that $n_0(x^1_0)=n_0$ and:
\[
g(x^2_0(l_0))=F_0\cup F_1\cup F_2\cup\bigcup_{i<\lh(w^2_1),i\neq l_1}(g(w^2_1(i))\cap A_{f(t_0^\frown (i))})\cup\bigcup_{i<l_2}(g(w^2_2(i))\cap A_{f(t_1^\frown (i))}).
\]
Then $E_0(x^2_0)=E_2(x^2_0,w^2_1\,^\frown 0_\infty,w^2_2 \,^\frown(g^{-1}(F_2))^\frown 0_\infty)\in B^1_1$ (since $g(0)=\emptyset$ and $(x^2_0,w^2_1\,^\frown 0_\infty,w^2_2 \,^\frown(g^{-1}(F_2))^\frown 0_\infty)\in [u^{2}_0]_{l_0}\times[u^{2}_1]_{l_1}\times[u^2_2]$). Find a basic set $B^2_0=\{A\subseteq 2^{<\omega}:\ G^2_0\subseteq A,\ H^2_0\cap A=\emptyset\}$, for finite $G^2_0$ and $H^2_0$, such that $E_0(x^2_0)\in B^2_0\subseteq U_{0,n_0,2}\cap B^1_1$. We can assume that $G^2_0$ contains:
\[
g(x^2_0(l_0))\cup \bigcup_{i<\lh(w^2_0),i\neq l_0}(g(w^2_0(i))\cap A_{f((i))}) \supseteq F_2.
\]
Find also any $\bar{x}^2_1=(x^2_{1,0},x^2_{1,1})\in [w^{2}_0]_{l_0}\times[w^{2}_1]_{l_1}$ such that $n_1(\bar{x}^1_1)=n_1$, $E_1(\bar{x}^2_1)\in B^2_0$ and:
\[
g(x^2_{1,1}(l_0))=F_1\cup F_2\cup\bigcup_{i<l_2}(g(w^2_2(i))\cap A_{f(t_1^\frown (i))}).
\]
There is $B^2_1=\{A\subseteq 2^{<\omega}:\ G^2_1\subseteq A,\ H^2_1\cap A=\emptyset\}$ such that $E_1(\bar{x}^2_1)\in B^2_1\subseteq U_{1,n_1,2}\cap B^2_0$. Finally, find $\bar{x}^2_2=(x^2_{2,0},x^2_{2,1},x^2_{2,2})\in [w^2_0]_{l_0}\times[w^2_1]_{l_1}\times[w^2_2]$ such that $n_2(\bar{x}^2_2)=n_2$, $E_2(\bar{x}^2_2)\in B^2_1$ and $g(x^2_{2,2}(l_2))=F_2$. There is $B^2_2=\{A\subseteq 2^{<\omega}:\ G^2_2\subseteq A,\ H^2_2\cap A=\emptyset\}$ such that $E_2(\bar{x}^2_2)\in B^2_2\subseteq U_{2,n_2,2}\cap B^2_1$. 

\medskip

By continuing this procedure we will obtain the required sequences $(l_k)_{k\in\omega}$, $(n_k)_{k\in\omega}$ and basic sets $B^m_k$, for $k\in\omega$ and $m\geq k$. This finishes the entire proof.

\end{proof}

%%%%%%%%%%%%%%%%%%%%%%%%%%%%%%%%%%%%%%%%%%%%%%%%%%
%%%%% REFERENCES
%%%%%%%%%%%%%%%%%%%%%%%%%%%%%%%%%%%%%%%%%%%%%%%%%%

\bibliographystyle{amsplain}
\bibliography{references}

\end{document}